\documentclass[12pt]{article} %[numbers 31.04.25-13.01.25-31.08.24]  %,a4paper
\usepackage{amsfonts,amssymb} 
\usepackage[cp1251]{inputenc}%\usepackage[english,russian]{babel}

\usepackage[margin=25mm, a4paper]{geometry}  % journal quality

\usepackage{tikz}

\def\n{\noindent}  \def\?#1{}
\def\IZ{{\mathbb{Z}}}  \def\IR{{\mathbb{R}}}
 \def\cA{{\cal A}}   \def\cM{{\cal M}}
        
      \def\la{\lambda}

   \def\phi{\varphi}  
\def\ep{\varepsilon}   \def\bdelete#1{}
    \def\v#1{\vec{#1}}

\def\mod1{\,({\rm mod\ } 1)\,}

\def\beq#1#2{\begin{equation} \label{#1} #2 \end{equation}}
\def\bea#1{\begin{eqnarray*} #1 \end{eqnarray*}} \def\a{\!\!\!&\!\!\!\!&}
\def\beaq#1#2{\label{#1} \begin{eqnarray} #2 \end{eqnarray}}
\def\toas#1{\stackrel{#1}{\longrightarrow}}

\def\proof{\smallskip \noindent {\bf Proof. \ }}       %start of proof
\def\blanksquare{\,\,\,$\sqcup\!\!\!\!\sqcap$}         %blank  square
\def\qed{\hfill\blanksquare\linebreak\smallskip\par}   %end of proof
\def\thname{Theorem}  \def\lmname{Lemma}    \def\prname{Proposition}
\def\dfname{Definition}  \def\crname{Corollary}  \def\rmname{Remark}
\def\exname{Example}   \def\conname{Conjecture}

\newtheorem{theorem}{\thname}[section]   %Numbering: Theorem--Other section
\newtheorem{lemma}{\lmname}[section]     %{lemma}[theorem]{Lemma}   subsection
\newtheorem{corollary}[lemma]{\crname}   %lemma

\newtheorem{conjecture}[lemma]{\conname} %lemma

\newtheorem{dftn}{\dfname}[section]
\newenvironment{definition}{\begin{dftn}\rm}{\end{dftn}} %section
\def\bdef#1{\begin{definition} #1 \end{definition}}
\newtheorem{rmrk}[lemma]{\rmname}
\newenvironment{remark}{\begin{rmrk}\rm}{\end{rmrk}}     %lemma

\makeatletter\def\fps@figure{htbp}\makeatother %figure pos: tbp - standard

\begin{document}

\title{Are prime numbers and quadratic residues random?}
\author{Michael Blank\thanks{
        Higher School of Modern Mathematics MIPT, 
        1 Klimentovskiy per., Moscow, Russia;}
        \thanks{National Research University ``Higher School of Economics'';
        e-mail: mlblank@gmail.com}
       }
\date{March 30, 2025} % \today}
\maketitle

\begin{abstract} 
Appeals to randomness in various number-theoretic constructions appear regularly 
in modern scientific publications. 
Such famous names as V.I. Arnold, M. Katz, Ya.G. Sinai, and T. Tao are just a few examples. 
Unfortunately, all of these approaches rely on various, although often very non-trivial and 
elegant, heuristics. A new analytical approach is proposed to address the issue of 
randomness/complexity of an individual deterministic sequence. 
This approach demonstrates the expected high complexity of quadratic residues and 
the unexpectedly low complexity in the case of prime numbers. 
Technically, our approach is based on a new construction of the dynamical 
entropy of a single trajectory, which measures its complexity, in contrast 
to classical Kolmogorov-Sinai and topological entropies, which measure the 
complexity of the entire dynamical system.
\end{abstract}

{\small\n
2020 Mathematics Subject Classification. Primary: 37A44; Secondary: 37A35, 11N05, 11K65.\\
Key words: entropy, ergodic theory, randomness, complexity, prime number, quadratic residue.
}
%http://aps.ecnu.edu.cn/UserFiles/File/MSC2020-Mathematica l% 37Bxx Topological dynamics
% 11N05 Distribution of primes;  11A15 Power residues; 11K65 Arithmetic functions in prob. number theory
% 37A35 Entropy and other invariants, isomorphism, classification in ergodic theory
% 37A44 Relations between ergodic theory and number theory;  81P17 Quantum entropies;
% 28D20 Entropy and other invariants;   37M25 Computational methods for ergodic theory
% 94A17 Measures of information, entropy

%{Comments to ArXiv:} The normalization of local entropy has been changed and a number of results and estimates have been improved.
%%%%%%%%%%%%%%%%%%%%%%%%%%%%%%%%

\section{Introduction}\label{s:intro}
Due to the obvious complexity of various number-theoretic constructions and 
the variety of patterns of numbers, both specialists in this field and other 
mathematicians interested in similar problems create purely random models 
to describe them. Let us mention Cramer's model (see, for example, \cite{Cr,Kac}), 
the substitution of random numbers in a series that defines the zeta function 
(see, for example, \cite{Tao}), the analysis of hidden periodicities in emerging 
sequences and geometric properties of the Poisson process (see \cite{Ar}), 
or the construction of a natural invariant measure concentrated on the set 
of square-free numbers (see \cite{SC}). 
These and other approaches will be discussed in Sections ~\ref{s:prime} and \ref{s:residues}.

The range of conclusions presented in the text is also interesting. 
From Arnold's denial of the randomness of quadratic residues, to Cramer's model 
which asserts the randomness of prime numbers, there seems to be a focus on 
highlighting specific properties, while disregarding others. 
In particular, Cramer's model relies heavily on the choice of the desired 
probability distribution, without taking into account the fact that this 
characteristic, although undoubtedly important, does not fully determine 
the random process.

It is natural to wonder how to distinguish a truly ``complex/random'' sequence from 
a ``simple/non-random'' one. To do this, consider a typical realization of a Bernoulli 
random process with equal probabilities $(1/2,1/2)$ of zeros and ones as the best 
candidate of the first type. The opposite type of candidate is again the Bernoulli process,
but with probabilities $(0,1)$. Purely deterministic analogues of these processes are 
trajectories of the doubling map ($x \to 2x \mod1$) and the halving map ($x \to x/2$).
As we will see, this is in many ways contrary to the ``algorithmic'' approach, in which 
the complexity of a sequence is interpreted through the complexity of the description 
of the process that generates the sequence (the latter is trivial in both last examples)

The purpose of this article is to provide a quantitative answer to the question of 
the complexity/randomness of a single purely deterministic sequence of points of 
the types discussed above. Such approaches are well known, but they all have 
serious drawbacks, as they are either non-constructive or completely arbitrary. 
Perhaps the most elegant among them was proposed by A. N. Kolmogorov. 
His main idea was to regard the sequence under study as the trajectory of 
a dynamical system and to reduce the complexity of the sequence to the analysis 
of the ``simplest'' dynamical system that produces it.

We will follow the same idea, but keeping in mind that different trajectories of 
the same dynamical system can exhibit qualitatively different behavior. 
We will introduce a new concept of ``local'' dynamical entropy $h_{loc}$. 
Unlike the Kolmogorov-Sinai metric entropy, it is independent of the choice of 
invariant measure and allows us to study even dynamical systems that do not 
have such a measure (see, for example, \cite{Bl17} for a discussion of such systems).

In the following, we will further refine this concept in order to analyze the complexity 
and randomness of individual binary sequences representing the number-theoretic 
constructions mentioned above.

In order to distinguish between static and dynamic entropy-like characteristics,  
we use the uppercase letter $H$ in the static case and the lowercase letter $h$ 
in the dynamic case. 
The entropy-like characteristic of randomness of a sequence $\v{x}$,  
denoted by $h(\v{x})$, takes values in the set of nonnegative real numbers, $\IR_+$. 
The value of $h(\v{x})=0$ is interpreted as a non-random sequence $\v{x}$, 
while a value greater than $0$ is interpreted as random. 
To simplify notation, we will use the following symbols: $\lim^-$ for $\liminf$ and 
$\lim^+$ for $\limsup$. We will omit the signs when the lower and upper limits are equal.

The paper is organized as follows. In Section~\ref{s:info}, we recall the classical 
definitions related to the concept of Shannon entropy of a discrete distribution and 
demonstrate that, from the point of view of small perturbations, this functional exhibits 
a number of unexpected features. We then move on to Section~\ref{s:dynEnt}, where 
we deal with various dynamical versions of the entropy concept, starting with the 
well-known Kolmogorov-Sinai metric entropy of a dynamical system. In order to take 
into account the complexity properties of individual trajectories, we introduce two 
completely new notions: local and information entropies for a single trajectory. 
In Section~\ref{s:Complexity}, we specify these new entropies for the case of 
a discrete phase space and introduce new measures of complexity/randomness for 
sequences of points in a finite alphabet. 
Finally, we compare these characterizations to some well-known approaches. 
The remaining part of the paper focuses on the application of these constructions. 
In Sections~\ref{s:prime} and ~\ref{s:residues}, we will prove the non-random 
nature of the set of prime numbers and the randomness of quadratic residues, respectively.

The author is grateful to L.Bassalygo, A.Kalmynin, S.Pirogov, A.Shen, M.Tsfasman, 
and A.Vershik for useful discussions of the issues raised in the article, 
as well as to the anonymous referee for valuable comments and suggestions.

\section{Information (Shannon) entropy}\label{s:info}

Let $(\Omega,\Sigma)$ be a measurable space with the Borel $\sigma$-algebra 
$\Sigma$ of measurable subsets. 

Throughout this section we assume that the space $\Omega$ is discrete, 
i.e., $\Omega:=\{\omega_1,\omega_2,\dots\}$, and $\Sigma:=2^\Omega$. 
Let $\v{p}:=(p_1,p_2,\dots)$ denote the distribution (non-necessarily probabilistic) 
on $\Omega$, i.e. $p_i\ge0~\forall i$, but $||\v{p}||:=\sum_ip_i=\v{p}(\Omega)$  
may differ from $1$. We also denote $\log(\cdot):=\log_2(\cdot)$, and 
$H(C):=-C\log C$ for any constant $C\ge0$, $H(0)=0$.

\bdef{The (Shannon) {\em entropy} of $\v{p}$ is defined as 
$H(\v{p}):=-\sum_i p_i\log p_i$.}

The following lemma collects several important (but little-known) observations that we will need later. 

\begin{lemma}\label{l:ineq} Let $r:=\#\Omega\le\infty$. Then  
\beq{e:H-max}{H(\v{p}) \le ||\v{p}||\log r + H(||\v{p}||) 
                                      = ||\v{p}||(\log r - \log||\v{p}||),}%
\beq{e:dist0}{H(\v{p}+\v{q}) = H(\v{p}) + |\log (e\inf_ip_i)| \cdot ||\v{q}|| 
                                            + o(||\v{q}||\cdot|\log(e\inf_ip_i)|)}%
\beaq{e:dist1}{ |H(\v{p}+\v{q}) - H(\v{p})| 
                \le \min\{H(\v{p}) + H(||\v{p}||) + H(||\v{q}||) + ||\v{q}||\log(r), \\  
                   ~ 2|\log(e\inf_ip_i)|\cdot ||\v{q}||\}  \nonumber}%
and this functional reaches its maximum value on the uniform distribution.
\end{lemma}
\proof The first claim follows from the well-known fact that in the case of finite 
$\Omega$, the uniform probabilistic distribution (denoted by $\v{p}^u$) maximizes entropy. 
Therefore, %
\bea{\log r = H(\v{p}^u) \a\ge H(\v{p}/||\v{p}||) 
             = -\sum_i \frac{p_i}{||\v{p}||} \log(\frac{p_i}{||\v{p}||}) \\
            \a= -\frac1{||\v{p}||} \sum_i p_i\log p_i + \frac1{||\v{p}||} \sum_i p_i \log(||\v{p}||) \\
            \a= \frac1{||\v{p}||} H(\v{p}) + \log(||\v{p}||) ,}%
which implies (\ref{e:H-max}).

The key point of the proof of the remaining inequalities is the observation that \\%
for each $t>0$ the first derivative of the function $t\log t$ is well defined, 
which for small enough $\ep$ gives 
the relation %
\beq{e:delent}{ H(t+\ep) - H(t) = -\ep\log et + o(\ep|\log et|) .} %
Applying (\ref{e:delent}) to the explicit formula for $H(\v{p})$ we get the result. \qed

\begin{corollary}\label{c:ineq} The functional $H: \ell_1(\Omega) \to\IR_+\cup \{\infty\}$  
is continuous on the entire space $\ell_1(\Omega)$ if and only if $\#\Omega<\infty$. 
\end{corollary}

This simple result is important from both theoretical and practical points of view; 
in particular, it demonstrates that Khinchin’s classical axiomatic entropy construction \cite{Kh} 
cannot be extended from finite space to an infinite one. In practice, this makes it possible
control the accuracy of entropy calculations when we know the distribution only approximately,
which will be very handy in Sections~\ref{s:prime},\ref{s:residues}. 

Surprisingly, Lemma~\ref{l:ineq} seems new. At least, I was not able to find 
results of this sort in numerous publications devoted to the concept of entropy.

In the sequel we will need the following estimates.

\begin{lemma}\label{l:pert}
Let $\v{p}^{(N)}:=\{p^{(N)}_i\}_{i=1}^N$ with $p^{(N)}_i\ge0$ and 
let $\ep_N:=||\v{p}^{(N)}||\in(0,1]$. %\sum_{i=1}^Np^{(N)}_i\le 1$. 
\begin{itemize}
\item[(a)] If $\exists C,\alpha\in\IR_+:~~ \ep_N  \le CN^{-\alpha} \quad \forall N\gg1$ 
               then $\lim\limits_{N\to\infty}H(\v{p}^{(N)}) = 0$. 
\item[(b)] If $\exists C\in\IR_+:~~ \ep_N  \le C/\log N \quad \forall N\gg1$ 
        then $\limsup\limits_{N\to\infty}H(\v{p}^{(N)}) \le C$, 
        and the upper bound can be achieved. 
\end{itemize}
\end{lemma}

\proof By inequality (\ref{e:H-max}), the entropy of a finite distribution $\v{q}:=\{q_i\}_{i=1}^N$ 
with a given sum $||\v{q}||>0$ reaches its maximum $||\v{q}|| \log (N/||\v{q}||)$ 
when the distribution is uniform (i.e. $q_i:=||\v{q}||/N~~\forall i$). 
Its minimum $H(||\v{q}||)$ is achieved for a single-point distribution.

In the case (a) for $N\gg1$ 
$$ H(\v{p}^{(N)}) \le \ep_N \log (N/\ep_N) \le CN^{-\alpha} \log(N^{\alpha+1}/C) \toas{N\to\infty}0, $$
which proves statement (a).

In the case (b) for $N\gg1$  %
\bea{ H(\v{p}^{(N)}) \le \ep_N \log (N/\ep_N) \a\le \frac{C}{\log N} \log(N\log N/C) \\
                             \a= C  (1 + \frac{\log(\log N)}{\log N}) - \frac{C\log C}{\log N}
                            \toas{N\to\infty}C. }
Again, by Lemma~\ref{l:ineq}(1), the upper limit is attained on a sequence of uniform distributions. 
This completes the proof.  \qed

\section{Dynamical entropy in ergodic theory}\label{s:dynEnt}

Let us give a brief account on classical approaches to the construction of entropy-like 
characteristics of a discrete time dynamical system, defined by a measurable map $f$ 
from a measurable space $(X,\Sigma,\mu)$ into itself. We start with the Kolmogorov-Sinai 
construction (for details, see, e.g., \cite{BPSJ}). 

\bdef{Given a pair of finite measurable partitions $\Delta, \Delta'$ of $(X,\Sigma,\mu)$ by their 
common {\em refinement} one means 
$\Delta \bigvee \Delta' := \{\Delta_{i}\cap \Delta'_{j}: ~ \mu(\Delta_{i}\cap\Delta'_{j})>0\}$.}

Let $\mu\in\cM_f$ (the set of all $f$-invariant measures). 
Making the refinement of $\{f^{-1}\Delta_{i}\}$ we again get a finite measurable partition 
which we denote by $f^{-1}\Delta$. 
The $n$-th refinement $\Delta^n$ of the partition $\Delta$ can be defined inductively 
$$\Delta^n:=\Delta^{n-1}\bigvee f^{-1}\Delta^{n-1}, ~~\Delta^0:=\Delta .$$

\bdef{The conditional {\em Kolmogorov-Sinai entropy} of a partition is defined as
$$h_\mu(f|\Delta):=\liminf\limits_{n\to\infty}\frac1n H_\mu(\Delta^n) 
                           =\lim\limits_{n\to\infty}\frac1n H_\mu(\Delta^n),$$ 
where $H_\mu(\Delta):=-\sum_{i}\mu(\Delta_{i})\ln\mu(\Delta_{i})$ is the entropy of the 
discrete distribution $\{\mu(\Delta_{i})\}$.}

\bdef{The {\em Kolmogorov-Sinai metric entropy} of the dynamical system $(f,X,\Sigma,\mu)$ 
is $h_\mu(f):=\sup_\Delta h_\mu(f|\Delta)$.}

Alternative approaches are known for continuous maps $f\in C^0(X,X)$, 
where $(X,\rho)$ is a compact metric space. 

\bdef{The $n$-th {\em Bowen metric} $\rho_n$ on $X$ is defined as   %
$\rho_n(u, v):=\max\left\{\rho\left(f^k(u), f^k(v)\right):\;k =0, \ldots,n-1 \right\}$.}

Let $B_\ep^n(x)$ be the open ball of radius $\ep$ in the metric $\rho_n$ around $x$. 

\bdef{The {\em Brin-Katok measure-theoretical entropy} of a measure $\mu\in \cM_f(X)$ 
at a point $u\in X$ is 
$h_\mu(f,u):=-\lim\limits_{\ep\to0}
              \lim\limits_{n\to\infty}\frac{1}{n}\log\mu(B_\ep^n(u))$.}%

Roughly speaking $h_\mu(f,u)$ measures the exponential
rate of decay of the measure of points that stay $\ep$-close to
the point $u$ under forward iterates of the map $f$.

\begin{theorem}\cite{BK}\label{t:BK} 
$h_\mu(f,u)$ is well defined for an ergodic measure $\mu$ and does not depend 
on $u$ for $\mu$-a.e $u\in X$.
\end{theorem}

A topological version (independent on the choice of the measure $\mu$) is  
available in the case of a continuous map $f$ (see, for example, \cite{Ka}). 
In this setting $\Delta:=\{\Delta_i\}_1^r$ is a {\em covering} of $X$ by open sets.
Define a transition matrix $M:=\{m_{ij}\}$, where $m_{ij}=1$ 
if $\Delta_i \cap f^{-1}\Delta_j \ne \emptyset$ and $=0$ otherwise.

Then on the Cantor set $X_M$ (the space of sequences with 
the alphabet $\cA:=\{1,2,\dots,r\}$ with the  transition matrix $M$)  
the left shift map $\sigma$ defines a symbolic dynamical system. 

$$ \v{x} = (x_1x_2\dots x_k\underbrace{x_{k+1}x_{k+1}\dots x_{k+n}}_{\v{w}} 
                 x_{k+n+1}\dots  x_N\dots) $$

Denoting by $A_\Delta^n$ the set of all admissible words $\v{w}$ of length $n$  
(i.e. different pieces of length $n$ of all trajectories of $(\sigma,X_M)$)  
and by $\#A_\Delta^n$ -- the number of such words, we set  
$$ h_{{\rm top}}(f|\Delta):=\liminf_{n\to\infty}\frac1n \log(\#A_\Delta^n)
                   =\lim_{n\to\infty}\frac1n \log(\#A_\Delta^n) . $$

\bdef{The {\em topological entropy} is defined as 
$h_{{\rm top}}(f):=\sup_\Delta h_{{\rm top}}(f|\Delta)$.}

Note that the construction of Kolmogorov-Sinai metric entropy 
(as well as of Brin-Katok entropy) are based on the choice of 
$f$-invariant measure $\mu$ (and depends on it), and the construction of 
topological entropy makes sense only for continuous maps. 
On the other hand, a general measurable dynamical system needs not 
to have even a single invariant measure (not speaking about the assumption 
on continuity). Discussion of dynamical systems having no invariant measures 
can be found, for example, in \cite{Bl17}.

To overcome these difficulties we propose yet another entropy-like constructions. 

Let $\Delta:=\{\Delta_i\}$ be a finite partition of $X$ by admissible\footnote{The choice of 
    admissible partitions for a general metric space is nontrivial and will be studied elsewhere. 
    A typical example is a partition into convex sets.} 
measurable sets.
We refer to the indices of $\Delta_i$ as an alphabet $\cA$, which needs not to be finite. 
We say that on a starting segment of length $N$ of a given trajectory 
$\v{x}:=(x_1,x_2,\dots)$ of our system there is a word $\v{w}:=(w_1,\dots,w_n)$ 
composed of the letters $w_i\in\cA$, if there is $i$ such that 
$$x_{i+j}\in \Delta_{w_j} ~~\forall j<n+1, i+j \le N .$$
Denote by $L(\v{x},\v{w},N)$ the number of occurrences of a word $\v{w}$ 
in the starting piece of length $N$ of the trajectory $\v{x}$, and let 
$\v{p}(\v{x},n,N)=(p_1,p_2,\dots)$ be a distribution (frequency) 
of all such words of length $n$. 

\bdef{By the conditional {\em local entropy} of the trajectory $\v{x}$ we mean 
$$ h_{{\rm loc}}^\pm(\v{x}|\Delta) 
            := \lim\nolimits_{n\to\infty}^{\pm} 
                \lim\nolimits_{N\to\infty}^{\pm} \frac1n H(\v{p}(\v{x},n,N)) .$$}
Here $\pm$ refers to the upper and lower limits, and 
$$ H(\v{p}(\v{x},n,N)) := -\sum_{i=1} p_i\log p_i $$ 
is the entropy of the distribution $\v{p}(\v{x},n,N)$. 

It is clear that this construction is something intermediate between metric and topological
entropy, but works for any measurable map and provides information about ``complexity'' 
of individual trajectories. The ``locality'' of $h_{{\rm loc}}$ is trajectory-wise, unlike 
the Brin-Katok entropy which is point-wise. 

In the case of a continuous mapping, a topological type approach to the same problem 
is known (see \cite{Br,ZL}), which with very minor modifications can be formulated for 
the measurable case discussed in this article.

Denote by $L(\v{x},n,N)$ the number of different words of length $n$ in
the starting piece of length $N$ of the trajectory $\v{x}$. 

\bdef{By the conditional {\em information entropy} of the trajectory $\v{x}$ we mean 
$$ h_{{\rm info}}^\pm(\v{x}|\Delta) 
   := \lim\nolimits_{n\to\infty}^{\pm} \lim\nolimits_{N\to\infty}^{\pm} \frac1n\log L(\v{x},n,N) .$$}

We postpone the analysis of the coincidences of the upper and lower limits in 
the above constructions in the general setting to a separate publication due 
to the following reasons. First, in all examples that we consider in the 
present work, the limits do indeed coincide. Second, the general sub-additive 
argument used in constructions of this kind do not work in the somewhat more 
general setting considered in Section~\ref{s:residues}.

Finally, we define the unconditional versions of the entropies under consideration as follows: %
\beq{e:B-ent}{ h_{{\rm loc}}^\pm(\v{x}) := \sup_\Delta h_{{\rm loc}}^\pm(\v{x}|\Delta), \quad 
    h_{{\rm info}}^\pm(\v{x}) := \sup_\Delta h_{{\rm info}}^\pm(\v{x}|\Delta) .}%

\bdef{We say that a sequence of points $\v{x}:=(x_1,x_2,\dots),~~x_i\in X$ is {\em typical} 
with respect to a probabilistic measure $\mu$, if 
$$\lim_{n\to\infty}\frac1n\sum_{i=1}^n 1_A(x_i)=\mu(A) \quad \forall A\in\Sigma .$$ }

In other words, the sequence $\v{x}$ is distributed according to the measure $\mu_{\v{x}}=\mu$. 

\begin{lemma}\label{l:per} For a periodic sequence $\v{x}$ the measure $\mu_{\v{x}}$ 
is well defined, and $h_{{\rm loc}}(\v{x}) = h_{{\rm info}}(\v{x}) = 0$. 
\end{lemma}
\proof The claim about the measure $\mu_{\v{x}}$ is obvious and we discuss only 
statements about the entropies. Consider a $\ell$-periodic sequence $\v{x}$. 
For each $n>\ell$ there are at most $\ell$ different sub-words of length $n$, each 
with the frequency $1/\ell+o(1/N)$. The local entropy 
of this sequence is equal to $\frac1n\log\ell\toas{n\to\infty}0$.
Similarly, $L(\v{x},n,N)\le\ell~~\forall (n,N)$, which implies the claim about the information entropy.  \qed 

An important question is the range of values of the functionals 
$h_{{\rm loc}}$ and $h_{{\rm info}}$, which is described in the following Lemma,
but its proof will be deferred to the next section.
\begin{lemma}\label{l:values-gen} 
$\forall \alpha,\beta\in[0,\infty]~~\exists \v{x}^\alpha, \v{x}^\beta$ such that 
$h_{{\rm loc}}(\v{x}^\alpha)=\alpha, ~h_{{\rm info}}(\v{x}^\beta)=\beta$.  
\end{lemma}

It is known (see, for example, \cite{Ka}) that under reasonably mild assumptions on the 
dynamical system $(f,X,\Sigma)$ we have $h_{{\rm top}}(f) = \sup_\mu h_\mu(f)$, 
where the supremum is taken over all ergodic $f$-invariant measures. Despite that the 
calculation of the local entropy does not depend on any free parameter (like an invariant 
measure in the case of the metric entropy), 
a similar connection can be established between $h_{{\rm info}}(\v{x})$ and $h_{{\rm loc}}(\v{x})$.
\begin{lemma}\label{l:ineq-entr} 
For any sequence $\v{x}$ we have
$$ h_{{\rm loc}}^\pm(\v{x}) \le h_{{\rm info}}^\pm(\v{x}) .$$
\end{lemma}

\proof For a given pair $n<N$ consider the distribution $\v{p}(\v{x},n,N)=(p_1,p_2,\dots)$. 
Since $L(\v{x},n,N)$ is the number of different words of length $n$ in
the starting piece of length $N$ of the trajectory $\v{x}$, we deduce that at most $L(\v{x},n,N)$ 
entries of $\v{p}(\v{x},n,N)$ are positive. Using that the entropy of a finite distribution 
reaches its maximum on the uniform distribution, we get
$$ H(\v{p}(\v{x},n,N)) := -\sum_{i=1} p_i\log p_i \le \log L(\v{x},n,N) ,$$ 
which implies the claim. \qed

If the measure $\mu_{\v{x}}$ is well defined, one can explore the connections between 
our newly defined entropy-like characteristics and more classical approaches.
However, this is beyond the scope of the present paper and will be studied in a separate article. 
Note also that our local entropy can be easily modified to work with bi-infinite trajectories 
of measurable semi-groups, which will also be studied elsewhere.

In the next section we apply the above construction to the study of the 
``complexity'' of sequences of numbers considered as trajectories of unknown 
dynamical systems, and then to some number-theoretic constructions.

\section{Complexity}\label{s:Complexity}
One of the main concepts of complexity theory, introduced by A.N.~Kolmogorov,  
was to reduce the question of the complexity of a given sequence of points $x$ to 
analysis of the complexity of a dynamical system admitting it 
(i.e. this sequence is the trajectory of such a system). Naturally, there are many 
dynamical systems admitting a given sequence, and Kolmogorov 
proposed to take into account the ``simplest'' among them. To implement this 
concept, he uses a universal Turing machine to describe any dynamical 
system living in a finite phase space, and the complexity of such 
machine is described by the minimum length of the program generating it \cite{Kol,ZL,Br}. 
The beauty of this approach comes at the cost of being completely non-constructive.
In practice, this type of complexity can only be calculated for some toy examples.

Apart from the non-constructiveness of this approach there are two other important issues. 
First, for a given sequence there might be no dynamical systems, admitting it as a trajectory. 
Second, different trajectories of the same dynamical system may demonstrate very 
different qualitative properties.

An alternative (let's call it algorithmic) concept (see \cite{ZiL,ZL,LNP,DLLM,KS1,Br}) 
is based on the idea to treat a given sequence as an unordered set of points 
and boils down to various dimension-like characteristics of this set. 
A serious disadvantage here is that important information about the order of items is lost. 
Moreover, unlike our approach, all variants of the algorithmic concept are constructed only 
in the case of a finite alphabet.

In what follows we try to make use of Kolmogorov's idea of the ``dynamical origin'' 
of the sequences under study. 

Let $\cA:=\{a_1,a_2,\dots,a_r\}$ be a finite collection of ``letters'' 
(to which we refer as an alphabet), 
equipped with the complete $\sigma$-algebra $\Sigma:=2^\cA$, and let 
$\v{x}:=\{a_{k_i}\}_{i\in\IZ_+}$ be a sequence composed of the letters from 
this alphabet.

For a map $f:\cA\to\cA$, admitting the sequence $\v{x}$ as a trajectory, we may 
apply the notions of the local and information entropies, defined by the 
relations (\ref{e:B-ent}). There are two important observations here. 
First, these notions do not depend on the choice of the map $f$ admitting $\v{x}$.  
Second, in the present setting the supremum over all admissible partitions  
can be easily calculated due to the non-negativity of conditional entropies 
and the existence of the most fine partition, which coincides with the 
partition into points. Thus we get %
\beq{e:Hloc}{ h_{{\rm loc}}^\pm(\v{x}) 
           := \lim\nolimits_{n\to\infty}^{\pm} \lim\nolimits_{N\to\infty}^{\pm} \frac1n H(\v{p}(\v{x},n,N)) .}%
\beq{e:Hinfo}{ h_{{\rm info}}^\pm(\v{x}) 
   := \lim\nolimits_{n\to\infty}^{\pm} \lim\nolimits_{N\to\infty}^{\pm} \frac1n\log L(\v{x},n,N) .}%

It is worth noting that a functional qualitatively similar to (\ref{e:Hloc}) is known in the literature 
on information theory under the name ``finite-state dimension'' (see \cite{DLLM,LNP,KS1} and 
further references therein). It was introduced in 2004 as a finite-state version of classical 
Hausdorff dimension, and it measures the lower asymptotic density of information in an 
infinite sequence over a finite alphabet. Another approach similar to (\ref{e:Hinfo}) was introduced 
in \cite{ZiL} under the name ``compression ratio''. 
Neither infinite and necessarily probabilistic distributions nor connections to any version of 
dynamical entropy has been discussed in this context.

\begin{lemma}\label{l:values} 
$\forall \alpha,\beta\in[0,1]~~\exists \v{x}^\alpha, \v{x}^\beta$ such that 
$h_{{\rm loc}}(\v{x}^\alpha)=\alpha, ~h_{{\rm info}}(\v{x}^\beta)=\beta$.  
\end{lemma}
\proof For a given $q\in[0,1]$ consider the Bernoulli random process ${\rm Ber}(q)$, 
that is a sequence of i.i.d. binary random variables $\{\xi_n\}_{n\in\IZ_+}$ with 
Prob$(\xi_n=1)=q$. Denoting by $\v{p}({\rm Ber}(q),n)$ the distribution of 
binary sub-words of length $n$ of ${\rm Ber}(q)$, we get (see \cite{Bi,Gr})%
\beq{e:Ber}{ \frac1n H(\v{p}({\rm Ber}(q),n)) 
                     = -q\log q - (1 - q)\log(1 - q) .} %
Moreover, by Shannon-McMillan-Breiman Theorem (see \cite{Bi}) for almost every 
realization $\v{x}^q$ of ${\rm Ber}(q)$ we have 
$$ \lim_{N\to\infty} \frac1n H(\v{p}(\v{x}^q,n,N)) 
     = \frac1n H(\v{p}({\rm Ber}(q),n)) ,$$
where $\v{p}(\v{x}^q,n,N)$ is the distribution of sub-words of length $n$ in the 
starting piece of length $N$ of $\v{x}^q$. The observation that the right hand side 
of (\ref{e:Ber}) depends on $q$ continuously and takes values in $[0,1]$ 
proves the first claim.

To prove the second claim, note that  $h_{{\rm info}}({\rm Ber}(q),n)\equiv\log2=1$ 
a.s. $\forall q\in(0,1)$. Consider instead of the Bernoulli process 
a topological Markov chain - collection of sequences consisted of letters from 
a finite alphabet $\cA$ defined by a binary transition matrix. 
Choosing elements of this matrix, one easily controls the 
number of different admissible sub-words of given length.\qed

The claims of Lemma~\ref{l:values} can be easily generalized for the case of the alphabet with 
an arbitrary (but finite) number of elements. 

\begin{lemma}\label{l:values-gen-r} 
Let $r:\#\cA<\infty$. Then 
$\forall \alpha,\beta\in[0,\log r]~~\exists \v{x}^\alpha, \v{x}^\beta$ such that 
$h_{{\rm loc}}(\v{x}^\alpha)=\alpha, ~h_{{\rm info}}(\v{x}^\beta)=\beta$.  
\end{lemma}
\proof Observe that for any positive integer $n$ the total number of 
all different words of length $n$ is equal to $r^n$. Therefore for any sequence $\v{x}$ 
with this alphabet by (\ref{e:H-max}) we get
$$ H(\v{p}(\v{x},n,N)) \le \log r^n = n\log r .$$
Using exactly the same arguments 
as in the proof of Lemma~\ref{l:values} we get that each value from the segment 
$[0,\log r]$ is admissible for our entropies. \qed

\n{\bf Proof} of Lemma~\ref{l:values-gen}. From the previous result we see that 
choosing an alphabet $\cA$ wiith $r$ symbols we can construct a sequence having 
entropies in the range $[0,\log r]$. Since $r$ is arbitrary this proves the claim. \qed

In the sequel we will pay a special attention to binary sequences $\v{b}$ with 
only 0 and 1 entries. Let us discuss their statistical properties in some detail. 

Denote by $M(\v{b},\v{w},N)$ the number of occurrences of the word $\v{w}$ among 
the first $N$ letters of $\v{b}$. 
By a zero word we will mean any locally maximal word, consisting of zeros only, 
and $\v1$ stands for the word consisting of a single letter $1$.

\begin{lemma}\label{l:rare}
Let $Q(\v{b},n,N)$ be the frequency of zero sub-words of length $n<N$ in the first $N$ letters 
of the sequence $\v{b}$. Then %
\beq{e:frec}{Q(\v{b},n,N) \ge 1 - \frac{M(\v{b},\v1,N)}{N} \frac{n}{1-n/N} .}%
\end{lemma}
\proof The binary sequence of length $N$, having the smallest number of zero sub-words 
of length $n$ can be realized as follows: $0\dots01~0\dots01~\dots~0\dots01~0\dots0$. 
Each of $M(\v{b},\v1,N)$ blocks $0\dots01$ consists of $(n-1)$ zeros and $1$ in the end. 
Thus the number of zero sub-words of length $n$ is equal to $N - nM(\v{b},\v1,N) - n$, 
while the total number of sub-words of length $n$ is $N-n$. Therefore the frequency 
$$ Q(\v{b},n,N) = \frac{N - nM(\v{b},\v1,N) - n}{N-n} = 1 - \frac{nM(\v{b},\v1,N)}{N-n} 
                        = 1 - \frac{M(\v{b},\v1,N)}{N} \frac{n}{1-n/N} .  $$ %\hskip1cm{\hbox{\blanksquare}}
Lemma is proven. \qed

\begin{lemma}\label{l:binary}
\begin{itemize}
\item[(a)] If $M(\v{b},\v1,N) \le CN^{1-\alpha}, ~\alpha\in(0,1)$, then 
        $Q(\v{b},n,N)\ge 1-2CnN^{-\alpha}$ and $h_{{\rm loc}}(\v{b})=0$. 
\item[(b)] If $M(\v{b},\v1,N) \le C\frac{N}{\log N}~\forall N\gg1$, then 
        $Q(\v{b},n,N)\ge 1-C\frac{n}{\log N}$ and $h_{{\rm loc}}(\v{b})\le C$. 
\end{itemize}
\end{lemma}
\proof Both claims follow from the direct application of Lemmas~\ref{l:rare} and \ref{l:pert}(a,b) 
respectively. \qed

\begin{remark} 
The claim of Lemma~\ref{l:binary}(1) can be misinterpreted as meaning that the zero 
density of ones in a binary sequence $\v{b}$ implies that $h_{{\rm loc}}(\v{b})=0$. 
To demonstrate that this is not the case, observe that the assumption in Lemma~\ref{l:binary}(2) 
allows the zero density of ones in $\v{b}$, but demonstrates that the local entropy 
in that case might be strictly positive. On the other hand, the high density of ones 
may lead to the zero local entropy as well (consider the sequence consisting of ones only). 
\end{remark}

Consider yet another interesting example. Let $\v{b}^{{\rm nat}}$ be the binary 
sequence obtained by concatenating the binary representations of all natural numbers. 
This sequence was introduced by D.G.~Champernowne \cite{Champ}. It is of interest, 
that the number whose fractional part coincides with the sequence 
$\v{b}^{{\rm nat}}$ is transcendental\footnote{Not the root of a non-zero polynomial 
     of finite degree with rational coefficients.} 
(see \cite{Mah}). 

\begin{lemma} $h_{{\rm loc}}(\v{b}^{{\rm nat}}) = h_{{\rm info}}(\v{b}^{{\rm nat}})=1$.
\end{lemma}
\proof It is known (see, for example, \cite{NaSh}) that the sequence $\v{b}^{{\rm nat}}$, 
considered as a binary number, is normal in any base\footnote{All its digits of the number 
represented in the given base follow the uniform distribution.}.
Therefore, the distribution of all binary sub-words of $\v{b}^{{\rm nat}}$ of the same length 
is uniform, which proves the first claim. The second claim follows from the observation that 
all possible finite binary words are present in this sequence. \qed

\section{Spatial distribution of prime numbers}\label{s:prime}

In publications on number theory (see, for example, \cite{Tao,Kac,Sour1,Sour2,Sour3,BFT}) 
we often read that random models provide heuristic support for various conjectures and 
that prime numbers are believed to behave pseudo-randomly in many ways and do not 
follow any simple pattern.
An important example of a purported statement about the pseudo-randomness of primes 
(known as the Cramer's model) is the Hardy-Littlewood conjecture for k-tuples, namely that 
the number of occurrences of different patterns in primes can be approximated by treating 
them as a sequence of random numbers generated by independently counting each $k \in\IZ_+$ 
as ``prime'' with probability $1/\log k$. For a detailed discussion of this and a number 
of other examples of this kind, see \cite{Sour3}. 

Consider the sequence of prime numbers $\v\pi:=(1, 2, 3, 5,7,11,13,\dots)$ and match it to the 
binary sequence $\v{b}^{{\rm prime}}:=\{b_i\}_{i\in\IZ_+}$, 
such that $b_{\pi_i}=1~\forall i\in\IZ_+$ and $b_{j}=0~\forall j\not\in\v\pi$.

\begin{theorem}\label{t:primes}  {\begin{itemize}
\item[(a)] $h_{{\rm loc}}(\v{b}^{{\rm prime}})
               \le h_{{\rm info}}(\v{b}^{{\rm prime}})\le \log((1+\sqrt5)/2) \approx 0.69424191363$, 
\item[(b)] under the validity of Hardy-Littlewood Conjecture (see Conjecture~\ref{c:tuples} below) 
               $h_{{\rm info}}(\v{b}^{{\rm prime}})=\log((1+\sqrt5)/2)$.
\end{itemize}}
%$h_{{\rm loc}}(\v{b}^{{\rm prime}})=0$. 
\end{theorem}

\begin{remark}
The distribution of finite patterns in $\v{b}^{{\rm prime}}$ is quite uneven, 
while the calculation of the upper bound for $h_{{\rm loc}}(\v{b}^{{\rm prime}})$ 
is based on the uniform one. 
Therefore, we expect that the true value of local entropy is much smaller and may 
even be zero. Moreover, below we will show that the local entropy of the Cramer's model 
of prime numbers is indeed zero.
\end{remark}

To prove Theorem~\ref{t:primes}  we need to discuss connections of some known 
statistics of primes to similar statistics of all finite binary words.

\begin{conjecture}\cite{IK}\label{c:tuples} 
Let $\v{a}:=(a_1,\dots,a_k)$ be distinct positive even integers which do not cover 
all residue classes to any prime modulus. Then the number of integers $0 < m \le N$ for which 
$m+a_1,\dots,m+a_k$ are all primes satisfies the asymptotic formula %
\beq{e:conj}{L_k(N,\v{a}) \approx C(k) N \log^{-k} N.  }
\end{conjecture}
%\cite{IK} (end of chapter 6, theorem 6.7, p.172). 

This conjecture is a close relative to the Hardy-Littlewood conjectures, but at present 
only partial results in this direction have been rigorously proven.

\bdef{A pair of consecutive prime numbers, separated by a single composite number 
are called {\em prime twins}.} 

Denote by $L_1(N)$ and $L_2(N)$ the number of primes and prime twins in $1,2,\dots,N$ correspondingly. 

\begin{theorem}\cite{Bar,IK}\label{t:twins} 
\begin{itemize}
\item $\frac{N}{\ln N - 2} < L_1(N) < \frac{N}{\ln N - 4}~~\forall N>54$. 
\item $L_2(N) < C \frac{N}{\ln^2N}$ for some $C<\infty$ and all $N\gg1$. 
\end{itemize}
\end{theorem}

Let us compare these statistics with similar information about general binary sequences. 
For this we need the following folklore result, well known in coding theory, but since we cannot 
give a precise reference, we will state it here and give a brief proof.

\begin{lemma}\label{l:twins} Let $Q(n,\v{w})$ be the number of all binary words of length $n$, 
avoiding the given word $\v{w}$. Then $Q(n,11) = Q(n-1,11) + Q(n-2,11)~~\forall n\ge3$. 
Therefore $Q(n,11)$ are Fibonacci numbers, satisfying the limit relations 
$z_n/z_{n+1}\to \la$, $z_{n+1}/(z_{n+1} + z_{n})\to \la$, where 
$\la^2+\la=1, ~~\la:=(1+\sqrt{5})/2\approx 1.61803398875$, and thus $z_n \approx \la^n$.
\end{lemma}
\proof Let's divide all words of length n without two 1's in a row into two groups: 
ending in 0 and ending in 1.
In the first case, any word of length $(n-1)$ without two 1's in a row can 
come before 0, and there are $Q(n-1,11)$ such words in total.
In the second case, 0 must come immediately before the last 1, and any of 
$Q(n-2,11)$ words of length $(n-2)$ without two 1's in a row can come before it.
From this we get $Q(n,11) = Q(n-1,11) + Q(n-2,11)$, as required. \qed

\begin{remark}
Similar (but more cumbersome) arguments can be used to find
statistics of words that avoid other combinations of zeros and ones
(such as the combination 101, which corresponds to prime twins).
We expect that analyzing these can improve our entropy estimates.
\end{remark}

\n{\bf Proof of Theorem~\ref{t:primes}.} 
Comparing the asymptotic number of prime numbers on large intervals $L_1(N)$ 
with the results of Lemma~\ref{l:binary}, we see that, despite the fact that prime numbers 
have zero density, they lie exactly on the boundary between sequences with 
zero local entropy and sequences with positivity entropy. 
Using part (b) of Lemma~\ref{l:binary}, we can obtain an upper bound that turns out to 
be quite large. 
Therefore, we will take a different approach, paying more attention to information entropy, 
which by Lemma~\ref{l:ineq-entr} cannot be smaller than the local entropy.

Note now that the sequence $\v{b}^{{\rm prime}}$ does not contain consecutive ones 
(except for the first three elements), which from the point of view of asymptotic relations 
do not play any role. 
This allows to make use of Lemma~\ref{l:twins} in order to estimate from above the total number 
of different sub-words of arbitrary length $n>2$ avoiding combinations 11. Namely,  we get %
\beq{e:prime}{ L(\v{b}^{{\rm prime}},n,N) \le C\left(\frac{1+\sqrt{5}}2\right)^n .}
This gives the desired upper bound for the information entropy:
$$ h_{{\rm info}}(\v{b}^{{\rm prime}}) 
   := \lim_{n\to\infty} \lim_{N\to\infty} \frac1n\log L(\v{b}^{{\rm prime}},n,N) 
   \le \log((1+\sqrt5)/2) .$$

The second statement about the exact value of information entropy follows from the assumption 
of the validity of the Hardy-Littlewood Conjecture, which states that there are all finite spatial 
combinations of prime numbers that avoid consecutive positions. In other words, only combinations 
of consecutive ones are excluded. Thus, the number of different sub-words is exactly equal 
(up to a multiplicative constant) to the right-hand side of the relation (\ref{e:prime}), 
which implies the result. \qed

\subsection{Inhomogeneous random Bernoulli process and Cramer's model}\label{s:cramer}

\bdef{{\em Bernoulli process} ${\rm Ber}(\v{q})$ with the vector-valued parameter 
$\v{q}:=\{q_k\}_{k=1}^\infty,~q_k\in[0,1]$ is a sequence of independent binary random variables 
$\{\xi_n\}_{n\in\IZ_+}$ with $P(\xi_n=1)=q_n$, where $P(\cdot)$ is the probability of an event.}

\begin{lemma}\label{l:Bernoulli} 
$H({\rm Ber}(\v{q}),n)=\sum_{k=1}^nH({\rm Ber}(q_k),1)$.
\end{lemma}
\proof Let $\v{b}^{(n)}$ be a binary word of length $n$. Then
$$ H(\v{p}({\rm Ber}(q),n)) := -\sum_{\v{b}^{(n)}} P(\v{b}^{(n)}) \log P(\v{b}^{(n)}) ,$$
where $P(\v{b}^{(n)})$ is the probability that the first $n$ letters of ${\rm Ber}(\v{q})$ coincide 
with $\v{b}^{(n)}$. On the other hand, setting $q'_{n+1}:=1-q_{n+1}$, 
due to the independence of the elements of the Bernoulli process, we obtain %
\bea{ H(\v{p}({\rm Ber}(\v{q}),n+1)) \a:= -\sum_{\v{b}^{(n)}} P(b^{(n+1)}) \log P(\v{b}^{(n+1)}) \\
  \a= -\sum_{\v{b}^{(n)}} P(\v{b}^{(n)})q_{n+1} \log (P(\v{b}^{(n)}) q_{n+1}) 
    - \sum_{\v{b}^{(n)}} P(\v{b}^{(n)})q'_{n+1} \log (P(\v{b}^{(n)}) q'_{n+1}) \\
  \a= -\sum_{\v{b}^{(n)}} P(\v{b}^{(n)})q_{n+1} \log P(\v{b}^{(n)}) 
        -\sum_{\v{b}^{(n)}} P(\v{b}^{(n)})q_{n+1} \log q_{n+1} \\
    \a~~~- \sum_{\v{b}^{(n)}} P(\v{b}^{(n)})q'_{n+1} \log P(\v{b}^{(n)}
       - \sum_{\v{b}^{(n)}} P(\v{b}^{(n)})q'_{n+1} \log q'_{n+1} \\
 \a= -(q_{n+1} + q'_{n+1}) \sum_{\v{b}^{(n)}} P(\v{b}^{(n)}) \log P(\v{b}^{(n)}) \\
     \a~~~   - (q_{n+1}\log q_{n+1} + q'_{n+1}\log q'_{n+1}) \sum_{\v{b}^{(n)}} P(\v{b}^{(n)}) \\ 
  \a= H(\v{p}({\rm Ber}(\v{q}),n)) + H(\v{p}({\rm Ber}(q_{n+1}),1)) , }%
since $q_{n+1} + q'_{n+1} =  \sum_{\v{b}^{(n)}} P(\v{b}^{(n)}) = 1$. \qed

In 1936, H.~Cramer \cite{Cr} introduced a probabilistic model of primes, where each natural number 
is selected for inclusion  with probability $1/\ln n$. From the point of view of the spatial 
distribution of these random numbers we are getting the inhomogeneous Bernoulli process 
${\rm Ber}(\v{q})$ with the vector-valued parameter $\v{q}:=\{q_k:=1/\ln k\}_{k=2}^\infty$.

\begin{theorem}\label{t:Cramer} 
For the Cramer's model $q_k:=\frac1{\ln k}$ we have $h_{{\rm loc}}({\rm Ber}(\v{q}))=0$.
\end{theorem}
\proof By definition,  
$$ h_{{\rm loc}}({\rm Ber}(\v{q})) := \lim_{n\to\infty} \frac1n H({\rm Ber}(\v{q}),n) .$$
On the other hand, by Lemma~\ref{l:Bernoulli} 
$$ \frac1n H({\rm Ber}(\v{q}),n) = \frac1n \sum_{k=2}^n H({\rm Ber}(q_k),1) 
            = -\frac1n\sum_{k=2}^n q_k\log q_k - \frac1n\sum_{k=2}^n q'_k\log q'_k .$$
We estimate the last two sums separately.
$$ S_n := -\sum_{k=2}^n q_k\log q_k = \sum_{k=2}^n \frac{\log\ln k}{\ln k}  
                                              < \sum_{k=2}^n \frac{k}{\ln k} .$$
Up to normalization, the last term is the mathematical expectation of the 
distribution $\{1/\ln k\}_{k=2}^n$. The function $1/\ln k$ decreases monotonically, 
so after normalization the last term cannot exceed $(n-1)/2$. 

To perform normalization, we need to calculate 
$$ R(n) := \sum_{k=2}^n \frac{1}{\ln k} 
             = \sum_{k=2}^n q_k = \sum_{k=2}^{\frac{n}{\ln^2n}}q_k 
             + \sum_{k=\frac{n}{\ln^2n}}^n q_k .$$
Clearly,
$$ \sum_{k=2}^{\frac{n}{\ln^2n}}q_k \le  \frac{n}{\ln^2n} .$$ 
On the other hand, 
$$\frac{n}{\ln n} < (n - \frac{n}{\ln^2n} +1) q_n < 
\sum_{k=\frac{n}{\ln^2n}}^n q_k \le \frac{n(1 - q_n)}{\ln n - \log\ln^2n} 
                                                      < \frac{n}{\ln n} + o(\frac{n}{\ln n}).$$
Therefore, $R(n) - \frac{n}{\ln n} = o(\frac{n}{\ln n})$, and $S_n \le \frac{n}{2\ln n}$. 

The 2nd sum boils down to 
$$ S'_n := -\sum_{k=2}^n q'_k\log q'_k 
             = -\sum_{k=2}^n (1 - \frac1{\ln k})  \log(1 - \frac1{\ln k}) $$
$$             = -\sum_{k=2}^n \log(1 - \frac1{\ln k}) 
             + \sum_{k=2}^n \frac1{\ln k} \log(1 - \frac1{\ln k}) $$
$$         \le C\sum_{k=2}^n \frac1{\ln k} \le C\frac{n}{\ln n} .$$

Finally, collecting above estimates, we get
$$ \frac1n H({\rm Ber}(\v{q}),n) = \frac1n C \frac{n}{\ln n} \toas{n\to\infty}0 .$$
\qed

In distinction to the homogeneous case one cannot use Shannon-McMillan-Breiman Theorem 
and we can claim that statistics of realizations $\frac1nH(\v{p}(\v{x}^{\v{q}},n,N))$ 
converge to $\frac1n H(\v{p}({\rm Ber}(\v{q}),n))$ only in probability.
Recall that $\v{p}(\v{x}^{\v{q}},n,N)$ is the distribution of sub-words of length $n$ in the 
starting piece of length $N$ of $\v{x}^{\v{q}}$.

\section{Spatial distribution of quadratic residues}\label{s:residues}

Patterns formed by quadratic residues and non-residues modulo a prime have been 
studied since the 19th century \cite{Ala} and still they continue to attract 
attention of contemporary mathematicians \cite{Ar,Co,MT,KTVZ} from various 
points of view. For a detailed historical overview of the concept of 
quadratic residues, we refer the reader to the monograph \cite{Wr}, 
and to the modern analysis from an algebraic-geometric point of view - to \cite{KTVZ}.

Probably V.I.~Arnold \cite{Ar} was the first to discuss these matters from 
the point of view of randomness, albeit at a heuristic level. 
Arnold's negative answer to this question stands in stark contrast 
to S.~Wright's \cite{Wr} positive answer, who used a completely 
different heuristic approach based on the Central Limit Theorem. 
S.~Wright argues also that the positive answer follows from earlier 
results due to P.~Kurlberg and Z.~Rudnick \cite{KR,Ku} about the 
distribution of spacings between quadratic residues. 

\bdef{An integer $k$ is called a {\em quadratic residue} modulo $q$ if it is congruent to 
a perfect square modulo $q$; i.e., if there exists an integer $\ell$ such that:
$\ell^2\equiv k~ ({\rm mod}~q)$. Otherwise, $k$ is called a quadratic non-residue modulo $q$.}

For $q = 19$, the 9 quadratic residues are $(1,4,5,6,7,9,11,16,17)$, 
while another 9 numbers $(2,3,8,10,12,13,14,15,18)$ are quadratic non-residues. 
A number of other examples with their analysis can be found in \cite{Wr}.

For an odd prime $q$ consider the finite sequence $(1, 2, \dots, q-1)$. 
Replacing each number in this sequence with $1$ if it is a quadratic residue modulo $q$, 
and $0$ otherwise, we get the binary word $\v{b}^{(q)}:=(b_1, \dots , b_{q-1})$.

In this section we will be interested in the ``randomness'' of the sequence of these 
binary words growing as $q\to\infty$ . Note that longer words here do not include 
shorter ones, and that none of the complexity-type concepts discussed in the literature 
capture situations of this kind. Therefore we need a new definition in terms of a scheme 
of series.

\bdef{Let $\v{B}:=(\v{b}^{(m)})_{m\in\IZ_+}$ be a sequence of binary words 
with $|\v{b}^{(m)}|\toas{m\to\infty}\infty$. By the {\em local entropy} of $\v{B}$
we mean $h_{{\rm loc}}(\v{B}):= \lim_{m\to\infty}h_{{\rm loc}}(\v{b}^{(m)})$, 
and by the {\em information entropy} 
$h_{{\rm info}}(\v{B}):= \lim_{m\to\infty}h_{{\rm info}}(\v{b}^{(m)})$.}

\begin{theorem}\label{t:residues} 
Let $\v{B}:=(\v{b}^{(q_m)})_{m\in\IZ_+}$ be a sequence of binary words representing 
quadratic residues, where $q_m$ is the $m$-th prime number. Then 
$h_{{\rm loc}}(\v{B}) = h_{{\rm info}}(\v{B})=1$. 
\end{theorem}

To prove this result we need some information about the statistics of quadratic residues.

Let $L(\v{w},\v{b})$ be the number of occurrences of the binary word $\v{w}$ in 
the finite binary sequence $\v{b}$. 
The following result on the asymptotic equidistribution of these quantities for 
quadratic residues was proved in \cite{Co} (see also a discussion in \cite{KTVZ}). 

\begin{theorem}\cite{Co}\label{t:co}
For each binary word $\v{w}$ of length $n:=|\v{w}|\le q$ we have 
$$ |L(\v{w},\v{b}^{(q)}) - q2^{-|\v{w}|}| < (|\v{w}|-1)\sqrt{q} + |\v{w}|/2 .$$
\end{theorem} 

\n{\bf Proof of Theorem~\ref{t:residues}}.
By Theorem~\ref{t:co} the distribution of sub-words of the same length in $\v{b}^{(q_m)}$ 
is asymptotically uniform with accuracy or order $1/\sqrt{q_m}$. 
Setting $N:=q_m, n:=|\v{w}|$, we see that the point-wise moduli of differences  
between the theoretical uniform distribution $\v{p}^{u}$ of words of length $n$ 
and the observed distribution $\v{p}^{o}$ cannot exceed $CN^{-1/2}$. 
Thus it seems that we are in a position to apply the result of Lemma~\ref{l:pert}.%{l:ineq}
Unfortunately this is not the case. Indeed, the $\ell_1$-norm of the difference of 
distributions is of order $NN^{-1/2}=N^{1/2}\toas{N\to\infty}\infty$.

Therefore, it is necessary to adjust the approach. Denote by 
$\{\ep_i:=p_i^{u} - p_i^{o}\}$ the point-wise differences of distributions. Then 
$$ |\ep_i| \le CN^{-1/2}, \quad \sum_i \ep_i = 0 .$$
The equality above follows from the fact the both distributions are probabilistic. 

Using the same arguments as in the proof of Lemma~\ref{l:ineq}, we get 
$$ |H(\v{p}^{u}) -  H(\v{p}^{o})| = \sum_i \ep_i |\log N| + o(1/N) = o(1/N) \toas{N\to\infty}0.$$

This proves that 
$$ h_{{\rm loc}}(\v{B}) = \lim_{N\to\infty} H(\v{p}^{u}) = 1 .$$ 
 
The last claim that $h_{{\rm info}}(\v{B})=1$ follows from the observation 
that by Theorem~\ref{t:co} all finite binary patterns have positive frequencies. \qed

%%%%%%%%%%%%%%%%%%%%%%%%%%%

%\newpage
  
\end{document}